\documentclass[11pt]{amsart}
\usepackage{amsmath}
\usepackage{amssymb}
\usepackage{amsthm}
\usepackage{psfig}
\newcommand{\ignore}[1]{\relax}

\newcommand{\C}{\mathbb C}
\newcommand{\R}{\mathbb R}
\newcommand{\Z}{\mathbb Z}

\newcommand{\ind}{\operatorname{ind}}

\newtheorem{thm}{Theorem}

\newtheorem{lem}[thm]{Lemma}
\newtheorem{cor}[thm]{Corollary}
\newtheorem{prop}[thm]{Proposition}
\newtheorem{add}[thm]{Addendum}
\theoremstyle{definition}
\newtheorem{defn}{Definition}
\newtheorem{que}{Question}
\newtheorem{exa}{Example}
\newtheorem{rem}{Remark}
\theoremstyle{remark}

\newtheorem{rmk}{Remark}

\newcommand{\tor}{(\C^*)^{n}}
\newcommand{\rtor}{(\R^*)^{n}}

\newcommand{\dd}{\partial}
\newcommand{\am}{\mathcal{A}}
\newcommand{\s}{\mathcal{S}}
\newcommand{\cp}{{\mathbb C}{\mathbb P}}
\newcommand{\rp}{{\mathbb R}{\mathbb P}}
\newcommand{\Log}{\operatorname{Log}}

\newcommand{\Vol}{\operatorname{Vol}}

\newcommand{\re}{\operatorname{Re}}

\newcommand{\Int}{\operatorname{Int}}
\renewcommand{\setminus}{\smallsetminus}

\newcommand{\Area}{\operatorname{Area}}

\begin{document}
\title{Amoebas of algebraic varieties}
\author{Grigory Mikhalkin}
\address{Dept of Math\\ Univ of Utah\\ Salt Lake City, UT 84112, USA}
\maketitle

The notion of amoebas for algebraic varieties was
introduced in 1994 by Gelfand, Kapranov and Zelevinski \cite{GKZ}.
Some traces of amoebas were appearing from time to time,
even before the formal introduction, as auxiliary tools in several problems
(see e.g. \cite{Be}). After 1994 amoebas have been
seen and studied in several areas of mathematics, from algebraic
geometry and topology to complex analysis and combinatorics.

In particular, amoebas provided a very powerful tool for studying
topology of algebraic varieties.
The purpose of this survey is to summarize our current state
of knowledge about amoebas and to outline their applications to
real algebraic geometry and adjacent areas. Most proofs are
omitted here. An expanded version of this survey is currently
under preparation jointly with Oleg Viro \cite{MV}.

\vspace{10pt}
\noindent {\em Acknowledgements.}
{This survey was prepared for the Real Algebraic and Analytic Geometry Congress,
June 11-15, 2001, Rennes, France. Most of the text was written during the stay of the author
at the Max-Planck-Institut, Bonn, Germany. The author thanks the Universit\'e de Rennes and
the Max-Planck-Institut for hospitality. The author is supported in part by the NSF.}

\section{Definition and basic properties of amoebas}
\subsection{Definitions}
Let $V\subset\tor$ be an algebraic variety. Recall that $\C^*=\C\setminus 0$ is the
group of complex numbers under multiplication.
Let $\Log:\tor\to\R^n$ be defined by $\Log(z_1,\dots,z_n)\to(\log|z_1|,\dots,\log|z_n|)$.
\begin{defn}[Gelfand-Kapranov-Zelevinski \cite{GKZ}]
The {\em amoeba} of $V$ is $\am=\Log(V)\subset\R^n$.
\end{defn}

\begin{figure}[h]
\label{amofline}
\centerline{
\psfig{figure=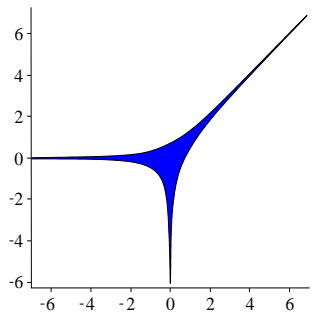,height=1.6in,width=1.6in}}
\caption{\cite{MV} The amoeba of the line $\{ x+y+1=0\} \subset (\C^*)^2$.}
\end{figure}

\begin{prop}[Gelfand-Kapranov-Zelevinski \cite{GKZ}]
The amoeba $\am\subset\R^n$ is a closed set with a non-empty complement.
\end{prop}

If $\C T\supset\tor$ is a closed $n$-dimensional toric variety and $\bar{V}\subset\C T$ is
a compactification of $V$ then we say that $\am$ is the amoeba of $\bar{V}$ (recall
that $\am$ is also the amoeba of $V=\bar{V}\cap\tor$).
Thus we can speak about amoebas of projective varieties once the coordinates
in $\cp^n$, or at least an action of $\tor$, is chosen.

If  $\C T$ is equipped with a $\tor$-invariant symplectic form then we can
also consider the corresponding moment map $\bar\mu:\C T\to\Delta$ (see \cite{At},\cite{GKZ}),
where $\Delta$ is the convex polyhedron associated to the toric variety $\C T$ with the
given symplectic form. The polyhedron $\Delta$ is a subset of $\R^n$ but it is well defined
only up to a translation.
In this case we can also define the {\em compactified amoeba} of $\bar{V}$.
\begin{defn}[Gelfand-Kapranov-Zelevinski \cite{GKZ}]
The {\em compactified amoeba} of $V$ is $\bar\am=\bar\mu(V)\subset\Delta$.
\end{defn}

\begin{figure}[h]
\centerline{\psfig{figure=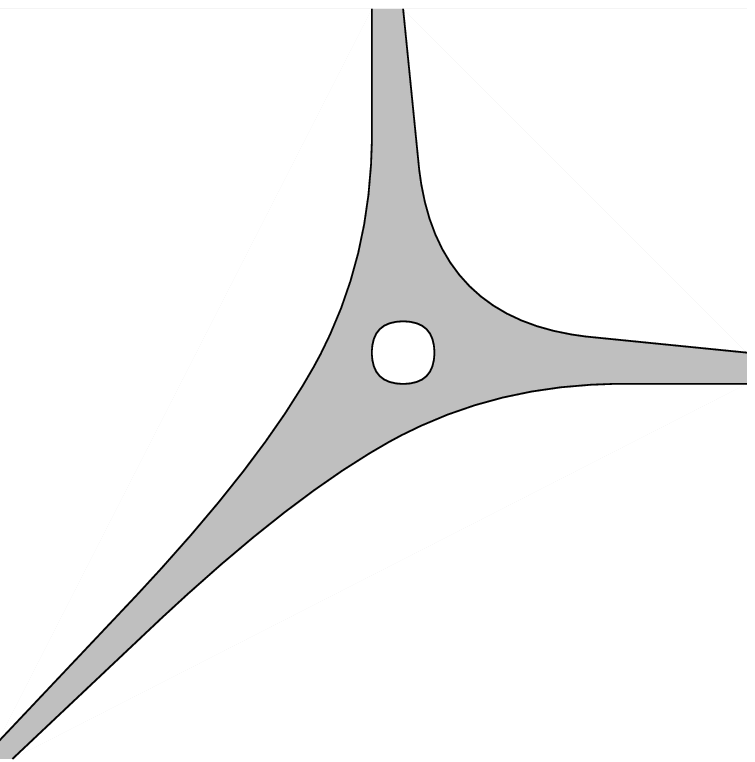,height=1in,width=1in}\hspace{1in}
\psfig{figure=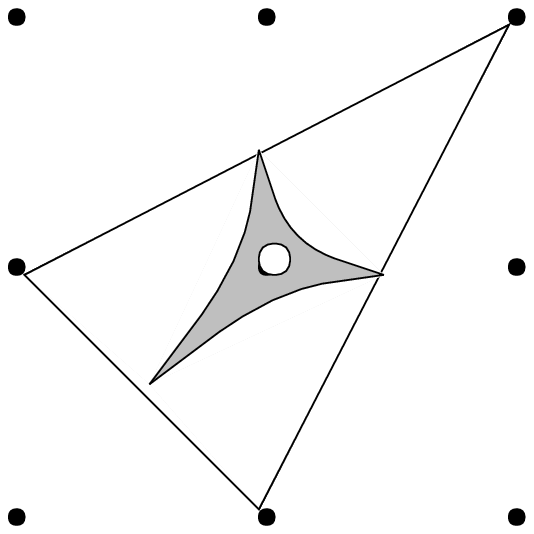,height=1in,width=1in}}
\caption{\cite{Mi} Amoeba $\am$ and compactified amoeba $\bar{\am}$.}
\end{figure}
\begin{rem}
\label{repar}
Maps $\bar\mu|_{\tor}$ and $\Log$ are submersions and have the same real $n$-tori as fibers.
Thus $\am$ is mapped diffeomorphically onto $\bar\am\cap\Int\Delta$ under a
reparametrization of $\R^n$ onto $\Int\Delta$.
\end{rem}

Using the compactified amoeba we can describe the behavior of $\am$
near infinity. Note that each face $\Delta'$ of $\Delta$ determines a
toric variety $\C T'\subset\C T$. Consider $\bar{V}'=\bar{V}\cap\C T'$.
Let $\bar\am'$ be the compactified amoeba of $\bar{V}'$.
\begin{prop}[Gelfand-Kapranov-Zelevinski \cite{GKZ}]
\label{natam}
We have
$\bar\am'=\bar\am\cap\Delta'$.
\end{prop}
This proposition can be used to describe the behavior of $\am\subset\R^n$ near infinity.

\subsection{Amoebas at infinity}
Consider a linear subspace $L\subset\R^n$ parallel to $\Delta'$ and
with $\dim L=\dim \Delta'$. Let $H\subset\R^n$ be a supporting hyperplane
for the convex polyhedron $\Delta$ at the face $\Delta'$, i.e. a hyperplane such
that $\Delta\cap H=\Delta'$.
Let $\stackrel{\to}{v}$ be an outwards normal vector to $H$.
Let $\am^{\Delta'}_t$, $t>0$, be the intersection of $L$ with the
result of translation of $\am$ by $-t\stackrel{\to}{v}$.
\begin{prop}
\label{asym}
The subsets $\am^{\Delta'}_t$ converge in the Hausdorff metric to $\am'$ when $t\to\infty$.
\end{prop}

\begin{figure}[h]
\label{3tent}
\centerline{
\psfig{figure=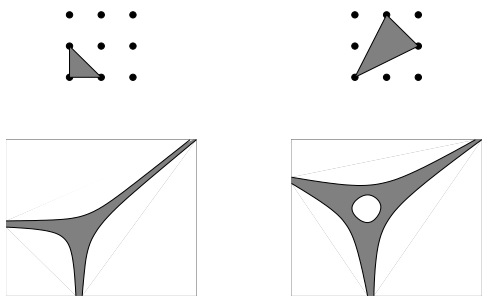,height=1.5in,width=2.5in}}
\caption{Amoebas together with their Newton polyhedra.}
\end{figure}

This proposition can be informally restated in the case $n=2$.
In this case $\Delta$ is a polygon and the amoeba $\am$ develops ``tentacles"
perpendicular to the sides of $\Delta$ (see Figure \ref{3tent}). The number
of tentacles perpendicular to a side of $\Delta$ is equal to the integer length
of this side, i.e. one plus the number of the lattice points in the interior of the side.

Note that we may assume (by passing to a different toric variety $\C T$ if needed)
that $V$ does not pass through the vertices of $\C T$, i.e. the fixed points of the $\tor$-action.
Thus we get the following corollary.
\begin{cor}
For a generic choice of the slope of a line $\ell$ in $\R^n$ the intersection $\am\cap\ell$ is compact.
\end{cor}

%
\subsection{Amoebas of hypersurfaces: concavity and topology of the complement}
\label{fpt}
Hypersurfaces case was treated by Forsberg, Passare and Tsikh in \cite{FPT}.
In this case $V$ is a zero set of a single polynomial $f(z)=\sum\limits_ja_jz^j,
a_j\in\C$. Here we use the multiindex notations $z=(z_1,\dots,z_n)$,
$j=(j_1,\dots,j_n)\in\Z^n$ and $z^j=z_1^{j_1}\dots z_n^{j_n}$.
Let
\begin{equation}
\label{NP}
\Delta=\text{Convex hull} \{j\ |\ a_j\neq 0\}\subset\R^n
\end{equation}
be the Newton polyhedron of $f$. 
\begin{thm}[Forsberg-Passare-Tsikh \cite{FPT}]
\label{thmfpt}
Each component of $\R^n\setminus\am$ is a convex domain in $\R^n$.
There exists a locally constant function
$$\ind:\R^n\setminus\am\to\Delta\cap\Z^n$$
which maps different components of the complement of $\am$ to
different lattice points of $\Delta$.
\end{thm}
\begin{cor}[Forsberg-Passare-Tsikh \cite{FPT}]
\label{chislokomponent}
The number of components of $\R^n\setminus\am$ is never greater
then the number of lattice points of $\Delta$.
\end{cor}
Theorem \ref{thmfpt} and Proposition \ref{asym} indicate the dependence
of the amoeba on the Newton polyhedron. 

The inequality of Corollary \ref{chislokomponent} is sharp.
This sharpness is a special case of Theorem \ref{komponenty} from
section \ref{paru}.
Also examples of amoebas with the maximal number
of the components of the complement are supplied by Theorem \ref{amlim}
from section \ref{secmas}.

The concavity of $\am$ is equivalent to concavity of its boundary.
The boundary $\dd\am$ is contained in the critical value locus of $\Log|_V$.
The following proposition also takes care of some interior branches of this
locus.
\begin{prop}[Mikhalkin \cite{Mi}]
\label{lc}
Let $D\subset\R^n$ be an open convex domain and $V'$ be a connected
component of $\Log^{-1}(D)\cap V$. Then $D\setminus\Log(V')$ is convex.
\end{prop}

\subsection{Amoebas in higher codimension: concavity}
The amoeba of a hypersurface is of full dimension in $\R^n$, $n>1$,
unless its Newton polyhedron $\Delta$ is contained in a line.
The boundary $\dd\am$ at its generic point is a smooth $(n-1)$-dimensional
submanifold. Its normal curvature form has no negative squares with
respect to the outwards normal (because
of convexity of components of $\R^n\setminus\am$). This property can be
generalized to the non-smooth points in the following way.
\begin{defn}
An open interval $D^1\subset L$, where $L$ is a straight line in $\R^n$,
is called a {\em supporting 1-cap} for $\am$ if
\begin{itemize}
\item $D^1\cap\am$ is non-empty and compact;
\item there exists a vector $\stackrel{\to}{v}\in\R^n$ such that the translation of $D^1$
by $\epsilon\hspace{-3pt}\stackrel{\to}{v}$ is disjoint from $\am$ for all sufficiently small $\epsilon>0$.
\end{itemize}
\end{defn}
The convexity of the components of $\R^n\setminus\am$ can be reformulated as stating
that {\em there are no 1-caps for $\am$}.

Similarly we may define higher-dimensional caps.
\begin{defn}
An open round disk $D^k\subset L$ of radius $\delta>0$ in a $k$-plane $L\subset\R^n$
is called a {\em supporting k-cap} for $\am$ if
\begin{itemize}
\item $D^k\cap\am$ is non-empty and compact;
\item there exists a vector $\stackrel{\to}{v}\in\R^n$ such that the translation of $D^k$
by $\epsilon\hspace{-3pt}\stackrel{\to}{v}$ is disjoint from $\am$ for all sufficiently small $\epsilon>0$.
\end{itemize}
\end{defn}

Consider now the general case, where $V\subset\tor$ is $l$-dimensional.
Let $k=n-l$ be the codimension of $V$.
The amoeba $\am$ is of full dimension in $\R^n$ if $2l\ge n$.
The boundary $\dd\am$ at its generic point is a smooth $(n-1)$-dimensional
submanifold. Its normal curvature form may not have more than
$k-1$ negative squares with respect to the outwards normal.
To see that note that a composition of $\Log|_V:V\to\R^n$ and
any linear projection $\R^n\to\R$ is a pluriharmonic function.

Note that this implies that there are no $k$-caps for $\am$
at its smooth points. It turns out that there are no $k$-caps
for $\am$ at the non-smooth points as well and also in the
case of $2l<n$ when $\am$ is $2l$-dimensional.
\begin{prop}[Local higher-dimensional concavity of $\am$]
\label{locconv}
If $V\subset\tor$ is of codimension $k$ then $\am$ does not
have supporting $k$-caps. 
\end{prop}

A global statement generalizing convexity of components  was recently
found by Andr\'e Henriques \cite{Henr}.
\begin{defn}[Henriques \cite{Henr}]
A subset $\am\subset\R^n$ is called {\em $k$-convex} if for any $k$-plane $L\subset\R^n$
the induced homomorphism $H_{k-1}(L\setminus\am)\to H_{k-1}(\R^n\setminus\am)$
is injective.
\end{defn}
\begin{thm}[Global higher-dimensional concavity of $\am$, cf. \cite{Henr}]
\label{globconv}
If $V\subset\tor$ is of codimension $k$ then $\am$ is $k$-convex.
\end{thm}
A proof of a weaker version of this statement is contained in \cite{Henr}.
Theorem \ref{globconv} can be deduced from its local version,
Proposition \ref{locconv}.

\subsection{Amoebas in higher codimension: topology of the complement}
Recall that in the hypersurface case each component of $\R^n\setminus\am$
is connected and that there are not more than $\#(\Delta\cap\Z^n)$ such
components. The correspondence between the components of the complement
and the lattice points of $\Delta$ can be viewed as a cohomology class
$\alpha\in H^0(\R^n\setminus\am;\Z^n)$ whose evaluation on a point in each
component of $\R^n\setminus\am$ is the corresponding lattice point.

Similarly, when $V$ is of codimension $k$ there exists a natural class (cf. \cite{R2})
$$\alpha\in H^{k-1}(\R^n\setminus\am;H^k(T^n)),$$
where $T^n$ is the real $n$-torus, the fiber of $\Log$, $H^k(T^n)=H^k(\tor)$.
The value of $\alpha$ on each $(k-1)$-cycle $C$ in $\R^n\setminus\am$ and $k$-cycle $C'$
in $T^n$ is the linking number in $\C^n\supset\tor$ of $C\times C'$ and the closure of $V$.

The cohomology class $\alpha$ corresponds to the linking with the fundamental class
of $V$. Consider now the linking with smaller-dimensional homology of $V$. 

Note that for an $l$-dimensional variety $V\subset\tor$ we have
 $H_j(V)=0$, $j>l$.
Similarly, $H^c_j(V)=0$, $j<l$, where $H^c$ stands for homology with closed support. 
The linking number in $\R^n$ composed with $\Log:\tor\to\R^n$ defines the
following pairing
$$H^c_l(V)\times H_{k-1}(\R^n\setminus\am)\to\Z.$$
Together with the Poincar\'e duality between $H^c_l(V)$ and $H_l(V)$
this pairing defines the homomorphism
$$\iota: H_{k-1}(\R^n\setminus\am)\to H_l(V).$$
\begin{que}
\label{sparivanie}
Is $\iota$ injective?
\end{que}
Recall that a subspace $L\subset H_l(V)$ is called {\em isotropic} if the restriction
of the intersection form to $L$ is trivial.
\begin{prop}
\label{isotrop}
The image $\iota(H_{k-1}(\R^n\setminus\am))$ is isotropic in $H_l(V)$.
\end{prop}
\begin{rmk}
A positive answer to Question \ref{sparivanie} together with Proposition  \ref{isotrop}
would produce an upper bound for the dimension of $H_{k-1}(\R^n\setminus\am)$.
\end{rmk}
One may also define similar linking forms for $H_j(\R^n\setminus\am)$, $j\neq k-1$
(if $j>k-1$ then we can use ordinary homology $H_{n-j-1}(V)$
instead of homology with closed support) .

The answer to Question \ref{sparivanie} is currently unknown even in the
case when $V\subset(\C^*)^2$ is a curve. In this case $V$ is a Riemann surface
and it is defined by a single polynomial.
Let $\Delta$ be the Newton polygon of $V$. 
The genus of $V$ is equal to the number of lattice points strictly inside $\Delta$
(see \cite{Kh}) while the number of punctures is equal to the number of lattice points
on the boundary of $\Delta$). Thus the dimension of a maximal isotropic subspace
of $H_1(V)$ is equal to $\#(\Delta\cap\Z^2)$ and Question \ref{sparivanie} agrees
with Corollary \ref{chislokomponent} for this case.

\section{Some analysis on amoebas}
\label{paru}
This section outlines the results obtained by Passare and Rullg{\aa}rd
in \cite{PR}, \cite{R1} and \cite{R2}.

We assume that $V\subset\tor$ is a hypersurface in this section.
Thus $V=\{f=0\}$ for a polynomial $f:\tor\to\C$ and we can consider
$\Delta\subset\R^n$, the Newton polyhedron of $V$ (see \ref{fpt}).

\subsection{The Ronkin function $N_f$}
Since $f$ is a holomorphic function,
$\log|f|:\tor\setminus V\to\R$ is a pluriharmonic function.
Furthermore, if we set $\log(0)=-\infty$ then
we have a plurisubharmonic function
$$\log|f|:\tor\to\R\cup\{-\infty\},$$
which is, obviously, strictly plurisubharmonic over $V$.
Recall that a function $F$ in a domain $\Omega\subset\C^n$ is called 
plurisubharmonic if its restriction to any complex line $L$ is subharmonic,
i.e. the value of $F$ at each point $z\in L$ is smaller or equal than the
average of the value of $F$ along a small circle in $L$ around $z$.

Let $N_f:\R^n\to\R$ be the push-forward of $\log|f|$ under the map $\Log:\tor\to\R^n$,
i.e.
$$N_f(x_1,\dots,x_n)=\frac{1}{(2\pi i)^n}\int\limits_{\Log^{-1}(x_1,\dots,x_n)}\log|f(z_1,\dots,z_n)|
\frac{dz_1}{z_1}\wedge\dots\wedge\frac{dz_n}{z_n},$$
cf. \cite{Ro}.
This function was called {\em the Ronkin function} in \cite{PR}.
It is easy to see that it takes real (finite) values even over $\am=\Log(V)$ where
the integral is singular.

\begin{prop}[Ronkin-Passare-Rullg{\aa}rd \cite{PR}, \cite{Ro}]
\label{ronkin}
The function $N_f:\R^n\to\R$ is convex. It is strictly convex
over $\am$ and linear over each component of $\R^n\setminus\am$.
\end{prop}
This follows from plurisubharmonicity of $\log|f|:\tor\to\R$, its strict
plurisubharmonicity over $V$ and its pluriharmonicity in $\tor\setminus V$.
Indeed the convexity of a function in a connected real domain is just
a real counterpart of plurisubharmonicity. A harmonic function of one real variable has
to be linear and thus a function of several real variables is real-plurisubharmonic
if and only if it is convex. Over each connected component of $\R^n\setminus\am$
the function is
linear as the push-forward of a pluriharmonic function.

\begin{rem}
Note that just the existence of a convex function $N_f$, which is
strictly convex over $\am$ and linear over components of $\R^n\setminus\am$,
implies that each component of $\R^n\setminus\am$ is convex.
\end{rem}

Thus the gradient $\nabla N_f:\R^n\to\R^n$ is constant over each component $E$ of $\R^n\setminus\am$.
Recall the classical Jensen's formula in complex analysis
$$\frac{1}{2\pi i}\int\limits_{|z|=e^x}\log|f(z)|\frac{dz}{z}=Nx+\log|f(0)|-\sum\limits_{k=1}^{N}\log|a_k|,$$
where $a_1,\dots,a_N$ are the zeroes of $f$ in $|z|<e^x$, if $f(0)\neq 0$ and $f(z)\neq 0$ if $|z|=e^x$.
This formula implies that $\nabla N_f(E)\in\Z^n\cap\Delta$.
\begin{prop}[Passare-Rullg{\aa}rd \cite{PR}]
We have
$$\Int\Delta\subset\nabla N_f(\R^n)\subset\Delta,$$
where $\Int\Delta$ is the interior of the Newton polyhedron.
\end{prop}
Recall that Theorem \ref{thmfpt} associates a lattice point
to each component of $\R^n\setminus\am$.
\begin{prop}[Passare-Rullg{\aa}rd \cite{PR}]
We have $$\nabla N_f(E)=\ind(E)$$
for each component $E$ of $\R^n\setminus\am$.
\end{prop}

\subsection{The spine of amoeba}

Passare and Rullg{\aa}rd \cite{PR} used $N_f$ to define {\em the spine} of amoeba.
Recall that $N_f$ is piecewise-linear on $\R^n\setminus\am$ and convex in $\R^n$.
Thus we may define a superscribed convex linear function $N^{\infty}_f$ by letting
$$N^{\infty}_f=\max_E N_E,$$
where $E$ runs over all components of $\R^n\setminus E$ and $N_E:\R^n\to\R$ is the linear
function obtained by extending $N_f|_E$ to $\R^n$ by linearity.

\begin{defn}[Passare-Rullg{\aa}rd \cite{PR}]
\label{spine}
The spine $\s$ of amoeba is the corner locus of $N^{\infty}_f$, i.e.
the set of points in $\R^n$ where $N^{\infty}_f$ is not locally linear.
\end{defn}
Note that $\s\subset\am$ and that $s$ is a piecewise-linear polyhedral complex.
The following theorem shows that $\s$ is indeed a spine of $\am$ in the
topological sense.

\begin{figure}[h]
\label{am-spine}
\centerline{
\psfig{figure=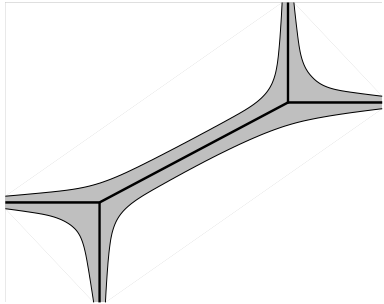,height=1.6in,width=2in}}
\caption{An amoeba and its spine.}
\end{figure}

\begin{thm}[Passare-Rullg{\aa}rd \cite{PR}, \cite{R2}]
The spine $\s$ is a strong deformational retract of the amoeba $\am$.
\end{thm}
Thus each component of $\R^n\setminus\s$ (i.e. each maximal open domain
where $N^{\infty}_f$ is linear) contains a unique component of $\R^n\setminus\am$.

\subsection{The spine $\s$ as a non-Archimedian amoeba}
Let $K$ be an arbitrary field with a norm. Let $K^*=K\setminus\{0\}$
and $V\subset(K^*)^n$ be an algebraic variety, i.e. the zero set of
a system of polynomial equations in $K$. The definition of amoeba
still makes sense in this setup: $\am_K=\Log(V)$, where $\Log:(K^*)^n\to\R^n$
is defined by $\Log(z_1,\dots,z_n)=(\log|z_1|,\dots,\log|z_n|)$.

Recall that a real-valued {\em valuation} of a field $K$ is a function $v:K^*\to\R$
such that $v(zw)=v(z)+v(w)$ and $-v(z+w)\ge\max\{-v(z),-v(w)\}$. A field with
such a valuation is called {\em non-Archimedian}. Note that $e^{-v}$ is
a (multiplicative) norm and $\Log:(K^*)^n\to\R^n$ is given by
$\Log(z_1,\dots,z_n)=(-v(z_1),\dots,-v(z_n))$.

Suppose that $K$ is algebraically closed and that $v:K^*\to\R$ is surjective.
\begin{exa}[cf. \cite{Ka}]
\label{pui}
Let $K$ be the real-power {\em Puiseux series} in $t$, i.e. the field whose elements
are formal power series $b(t)=\sum\limits_r b_rt^r$, where $b_r\in\C$ and the set of powers $r$
is bounded from below and is contained in a finite union of arithmetic progression.
This is an algebraically closed non-Archimedian field. The valuation of $b(t)$ is given by
the smallest power of $t$ which appears in the series.
\end{exa}

Unlike the complex case the amoeba of a hypersurface $V=\{f=0\}\subset(K^*)^n$
is completely determined by the norms of the coefficients of the defining polynomial $f$.
Let $f(z)=\sum\limits_j a_jz^j$, where $z\in(K^*)^n$ and $j\subset\Z^n$ is a multiindex.
A function $j\mapsto v(a_j))$ can be considered as a partially defined function on $\R^n$
(defined only on the finite set of lattice points $j$).
Its Legendre transform
$$N^K_f(x)=\max\limits_j \{jx-v(a_j)\},$$
where $jx$ is the scalar product of $j$ and $x$ in $\R^n$, is a convex piecewise-linear function $\R^n\to\R$.
\begin{thm}[Kapranov \cite{Ka}]
\label{thmka}
The amoeba $\am_K$ coincides with the corner locus of the piecewise-linear function $N^K_f$
(cf. Definition \ref{spine}). In particular, $\am_K$ is completely determined by the norms
of the coefficients of $f$.
\end{thm}

\subsection{Non-Archimedian amoebas as a counterpart of algebraic hypersurfaces}
\label{nA}
Subsets $\am_K\subset\R^n$ may be treated in a similar way we treat
algebraic hypersurfaces in $\tor$. Theorem \ref{thmka} ensures that
the choice of non-Archimedian field $K$ is irrelevant here
as long as $K$ is  algebraically closed and its valuation is onto $\R$.

Let us fix a Newton polyhedron $\Delta$. The space of all complex polynomials
which have $\Delta$ as its Newton polyhedron is $\C^N$, where $N=\#(\Delta\cap\Z^n)$.
Polynomials which are different by multiplication by a constant give the same
hypersurface. Thus the hypersurfaces space is $\cp^{N-1}$.

By Theorem \ref{thmka} the amoeba $\am_K$ is determined solely
by the valuations of the coefficients of the polynomial defining $V\subset (K^*)^n$.
A valuation on the monomials from $\Delta$ are functions
$\Delta\cap\Z^n\to\R$. They form the space $\R^N$.
Valuation functions which are different by adding a constant give
the same non-Archimedian amoebas. Thus non-Archimedian amoebas
are parametrized by $\R^{N-1}$.

\begin{rem}
In general, the space of non-Archimedian amoebas is not $\R^{N-1}$ but a quotient of $\R^{N-1}$.
If the valuation function $v:\Delta\cap\Z^n\to\R$ is not convex then we may vary
a little $v$ at some points keeping the Legendre transform the same.
\end{rem}

\begin{rem}[Non-Archimedian amoebas and Enumerative Geometry]
In a seminar talk in Paris, November 2000, Kontsevich noted a possibility of
using non-Archimedian amoebas in enumerative geometry.
As an example consider the problem of counting the number $n_d$ of
rational curves of degree $d$ in $\cp^2$ which
pass through $3d-1$ fixed generic points.
A generic complex  polynomial defines
a curve of genus $(d-1)(d-2)/2$. The polynomials defining rational
curves form a subset of codimension $(d-1)(d-2)/2$ and thus
the rational curves form a $(3d-1)$-dimensional space (the space of
curves has dimension one less than the dimension of the space of
corresponding polynomials).

\begin{figure}[h]
\label{na-cubic}
\centerline{
\psfig{figure=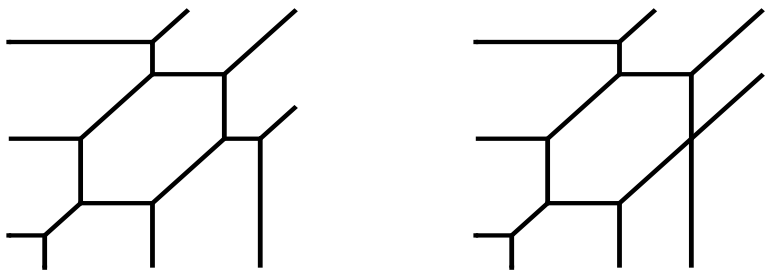,height=1.4in,width=4in}}
\caption{A smooth ``non-Archimedian cubic amoeba" and a rational ``non-Archimedian cubic amoeba".}
\end{figure}

Consider the space of non-Archimedian amoebas corresponding to curves of degree $d$ in ${\mathbb P}^2$.
This means that the Newton polygon $\Delta$ is the triangle with vertices $(0,0)$, $(d,0)$ and $(0,d)$.
It may be deduced from Theorem \ref{thmka} that a generic amoeba is a 3-valent graph which is
homotopy equivalent to a wedge of up to $(d-1)(d-2)/2$ circles.
Amoebas with $(d-1)(d-2)/2$ 4-valent vertices form a subset of codimension $(d-1)(d-2)/2$.
These amoebas play the role of rational curves.
Through generic $(3d-1)$ points in $\R^2$ there are $n_d$ of different amoebas of this kind.

As an exercise the reader may check that there is a unique non-Archimedian amoeba
of degree 1 through any 2 generic points in $\R^2$.
The 2 points are special for this problem if they belong to the same horizontal,
vertical or slope 1 line in $\R^2$. There is an infinite number of degree 1 non-Archimedian
amoebas through the 2 points if they are special.
\end{rem}

\subsection{Spine of amoebas and some functions on the space of complex polynomials}
Now we return to the study of the spine $\s\subset\am$ of a complex amoeba.
The spine $\s$ itself a certain amoeba over a non-Archimedian
field $K$. It does not matter what is the field $K$ as long as the corresponding hypersurface
over $K$ has the coefficients $a_j\in K$ with the correct valuations.
We can find these valuations from $N^{\infty}_f$ by taking its Legendre transform.
Since $N^{\infty}_f$ is obtained as a maximum of a finite number of linear function with
integer slopes its Legendre transform has a support on a convex lattice polyhedron $\Delta\subset\R^n$.
Let $c_{\alpha}\in\R$, $\alpha\in\Delta\cap\Z^n$ be the value of the Legendre transform of $N^{\infty}_f$ at $\alpha$.
To present $\s$ as a non-Archimedian amoeba we choose $a_j\in K$ such that $v(a_j)=c_{\alpha}$.

For each $\alpha\in\Delta\cap\Z^n$ let $U_{\alpha}$ be the space of all polynomials
whose Newton polyhedron is contained in $\Delta$ and whose amoeba contains a component
of the complement of index $\alpha$. The space of all polynomials whose Newton polyhedron
is contained in $\Delta$ is isomorphic to $\C^N$, where $N=\#(\Delta\cap\Z^n)$.
The subset $U_{\alpha}\subset\C^N$ is an open domain.
Note that $c_{\alpha}$ defines a real-valued function
on $U_{\alpha}$. This function was used by Rullg{\aa}rd \cite{R1}, \cite{R2} for
the study of geometry of $U_{\alpha}$.

\subsection{Geometry of $U_\alpha$}
Fix $\alpha\in\Delta\cap\Z^n$. Consider the following function in the space $\C^N$
of all polynomials $f$ whose Newton polyhedron is contained in $\Delta$
$$u_\alpha(f)=\inf\limits_{x\in\R^n}\frac{1}{(2\pi i)^n}
\int\limits_{\Log^{-1}(x)}\log|\frac{f(z)}{z^{\alpha}}|\frac{dz_1}{z_1}\wedge\dots\wedge\frac{dz_n}{z_n}, z\in\tor.$$
 Rullg{\aa}rd \cite{R1} observed that this function is plurisubharmonic in $\C^N$ while
pluriharmonic over $U_\alpha$. Indeed, over $U_\alpha$ there is a component
$E_\alpha\subset\R^n\setminus\am$ corresponding to $\alpha$ and $u_\alpha=\re\Phi_\alpha$,
where 
$$\Phi_\alpha=\frac{1}{(2\pi i)^n}\int\limits_{\Log^{-1}(x)}\log(\frac{f(z)}{z^{\alpha}})
\frac{dz_1}{z_1}\wedge\dots\wedge\frac{dz_n}{z_n}, x\in E_\alpha$$
is a $(\C/2\pi i\Z)$-valued holomorphic function. Note that over $\Log^{-1}(E_\alpha)$
we can choose a holomorphic branch of $\log(\frac{f(z)}{z^{\alpha}})$ and that
$\Phi_\alpha$ does not depend on the choice of $x\in E_\alpha$.
Therefore, $U_\alpha$ is pseudo-convex.

Note that $U_\alpha$ is invariant
under the natural $\C^*$-action in $\C^N$. Let ${\mathcal C}\subset\cp^{N-1}$ be the complement
of the image of $U_\alpha$ under the projection $\C^N\to\cp^{N-1}$.
\begin{thm}[Rullg{\aa}rd \cite{R1}]
For any line $L\subset\cp^{n-1}$ the set $L\cap {\mathcal C}$ is non-empty and connected.
\end{thm}

The next theorem describes how the sets $U_\alpha$ with different $\alpha\in\Delta\cap\Z^n$
intersect. It turns out that for any choice of subdivision $\Delta\cap\Z^n=A\cup B$ with $A\cap B=\emptyset$
the sets $\bigcup\limits_{\alpha\in A} U_\alpha$ and $\C^N\setminus\bigcup\limits_{\beta\in B} U_\beta$
intersect. A stronger statement was found by Rullg{\aa}rd.
Let $A, B\subset\Delta\cap\Z^n$ be disjoint sets. The set $A\cup B\subset\Delta\cap\Z^n$ defines
a subspace $\C^{\#(A\cup B)}\subset\C^N$ .
\begin{thm}[Rullg{\aa}rd \cite{R1}]
\label{komponenty}
For any $\#(A\cup B)$-dimensional space $L$ parallel to $\C^{\#(A\cup B)}$ the intersection
$L\cap\bigcup\limits_{\alpha\in A} U_\alpha\cap\C^N\setminus\bigcup\limits_{\beta\in B} U_\beta$
is non-empty.
\end{thm}

\subsection{The Monge-Amp\`ere measure and the symplectic volume}

\begin{defn}[Passare-Rullg{\aa}rd \cite{PR}]
\label{ma}
The Monge-Amp\`ere measure on $\am$ is the pull-back of the Lebesgue measure on $\Delta\subset\R^n$
under $\nabla N_f$.
\end{defn}
Indeed by Proposition \ref{ronkin} the Monge-Amp\`ere measure is well-defined. Furthermore,
we have the following proposition.

\begin{prop}[Passare-Rullg{\aa}rd \cite{PR}]
The Monge-Amp\`ere measure has its support on $\am$.
The total Monge-Amp\`ere measure of $\am$ is $\Vol\Delta$.
\end{prop}

By Definition \ref{ma} the Monge-Amp\`ere measure is given by the determinant of the Hessian of $N_f$.
By convexity of $N_f$ its Hessian $\operatorname{Hess} N_f$
is a non-negatively defined matrix-valued function.
The trace of $\operatorname{Hess} N_f$ is the Laplacian of $N_f$,
it gives another natural measure supported on $\am$.
Note that $\omega=\sum\limits_{k=1}^n\frac{dz_k}{z_k}\wedge\frac{d\bar z_k}{\bar z_k}$
is a symplectic form on $\tor$ invariant with respect to the group structure.
The restriction $\omega|_V$ is a symplectic form on $V$. Its $(n-1)$-th power
divided by $(n-1)!$ is
a volume form called {\em the symplectic volume} on the $(n-1)$-manifold $V$.
\begin{thm}[Passare-Rullg{\aa}rd \cite{PR}]
The measure on $\am$ defined by the Laplacian of $N_f$ coincides with the
push-forward of the symplectic volume on $V$, i.e. for any Borel set $A$
$$\int\limits_A\Delta N_f=\int\limits_{\Log^{-1}(A)\cap V}\omega^{n-1}.$$
\end{thm}

This theorem appears in \cite{PR} as a particular case of a computation for
{\em the mixed Monge-Amp\`ere operator}, the symmetric multilinear
operator associating a measure to $n$ functions $f_1,\dots,f_n$ (recall that
by our convention $n$ is the number of variables) and such
that its value on $f,\dots,f$ is the Monge-Amp\`ere measure from
Definition \ref{ma}. The total mixed Monge-Amp\`ere measure for $f_1,\dots,f_n$
is equal to the mixed volume of the Newton polyhedra of $f_1,\dots,f_n$
divided by $n!$.

Recall that this mixed volume divided by $n!$ appears in the Bernstein formula \cite{B}
which counts the number of common solutions of the system of equations $f_k=0$
(assuming that the corresponding hypersurfaces intersect transversely).
Passare and Rullg{\aa}rd found the following local analogue of the Bernstein formula
which also serves as a geometric interpretation of the mixed Monge-Amp\`ere measure.
Note that the complex torus $\tor$ acts on polynomials of $n$ variables. The value
of $t\in\tor$ on $f:\tor\to\C$ is the composition $f\circ t$ of the multiplication by $t$ followed by
application of $f$. In particular, the real torus $T^n=\Log^{-1}(0)\subset\tor$ acts on
polynomials of $n$ variables.
\begin{thm}[Passare-Rullg{\aa}rd \cite{PR}] 
\label{be}
The mixed Monge-Amp\`ere measure for $f_1,\dots,f_n$ of a Borel set $A\subset\R^n$
is equal to the average number of solutions of the system of equations $f_k\circ t_k=0$
in $\Log^{-1}(E)\subset\tor$, $t_k\in T^n$, $k=1,\dots,n$. 
\end{thm}
The number of solution of this system of equations does not depend on $t_k$
as long as the choice of $t_k$ is generic. Thus Theorem \ref{be} produces the Bernstein
formula when $E=\R^n$. 

\subsection{The area of a planar amoeba}
The computations of the previous subsection can be used to obtain an upper
bound on amoeba's area in the case when $V\subset(\C^*)^2$ is a curve.
With the help of Theorem \ref{be} Passare and Rullg{\aa}rd \cite{PR} showed
that in this case the Lebesgue measure on $\am$ is not greater than $\pi^2$
times the Monge-Amp\`ere measure. In particular we have the following theorem.
\begin{thm}[Passare-Rullg{\aa}rd \cite{PR}]
\label{PR}
If $V\subset(\C^*)^2$ is an algebraic curve then
$$\operatorname{Area}\am\le\pi^2\operatorname{Area}\Delta.$$
\end{thm}
This theorem is specific for the case $\am\subset\R^2$.
Non-degenerate higher-dimensional amoebas of hypersurfaces have infinite volume.
This follows from Proposition \ref{asym} since the area of the cross-section at infinity
must be separated from zero.

\section{Applications to real algebraic geometry}
\subsection{The first part of Hilbert's 16th problem}
\label{H16}
Most applications considered here are in the framework of Hilbert's 16th problem.
Consider the classical setup of its first part, see \cite{Hi}.
Let $\R\bar{V}\subset\rp^2$ be a smooth algebraic curve of degree $d$.
{\em What are the possible topological types of pairs $(\rp^2,\R\bar{V})$ for a given $d$?}

Since $\R \bar{V}$ is smooth it is homeomorphic to a disjoint union of circles.
All of these circles must be contractible in $\R P^2$ (such circles are called {\em the ovals})
if $d$ is even. If $d$ is odd then exactly one of these circles is non-contractible.
Therefore, the topological type of $(\rp^2,\R\bar{V})$ (also called {\em the topological
arrangement} of $\R\bar{V}$ in $\rp^2$) is determined by the number of components
of $\R\bar{V}$ together with the information on the mutual position of the ovals.

The possible number of components of $\R\bar{V}$ was determined by Harnack \cite{Ha}.
He proved that it cannot be greater than $\frac{(d-1)(d-2)}{2}+1$. Furthermore he proved
that for any number $$l\le\frac{(d-1)(d-2)}{2}+1$$ there exists a curve of degree $d$ with
exactly $l$ components as long as $l>0$ in the case of odd $d$ (recall that for odd $d$
we always have to have a non-contractible component).

Note that each oval separates $\rp^2$ into its {\em interior}, which is homeomorphic to a disk,
and its {\em exterior}, which is homeomorphic to a M\"obius band. If the interiors of the
ovals intersect then the ovals are called {\em nested}. Otherwise the ovals are called
{\em disjoint}. Hilbert's problem started from a question whether a curve of degree 6
which has 11 ovals (the maximal number according to Harnack) can have all of the ovals
disjoint. This question was answered negatively by Petrovsky \cite{P} who showed
that at least two ovals of a sextic must be nested if the total number of ovals is 11.
 
In general the number of topological  arrangements of curves of degree $d$
 grows exponentially with $d$. Even for small $d$ the number of the possible
types is enormous. Many powerful theorems restricting possible topological
arrangements were found for over 100 years of history of this problem,
see, in particular, \cite{P}, \cite{A}, \cite{R}. A powerful {\em patchworking}
construction technique \cite{V} counters these theorems.
The complete classifications is currently known for $d\le 7$, see \cite{V}.

The most restricted turn out to be curves with the maximal numbers of
components, i.e. with $l=\frac{(d-1)(d-2)}{2}+1$. Such curves were called
{\em M-curves} by Petrovsky. However, even for M-curves, the number of topological
arrangements grows exponentially with $d$. 

The situation becomes different if we consider $\rp^2$ as a toric surface,
i.e. as a compactification of $(\R^*)^2$. Recall that $\rp^2\setminus (\R^*)^2$
consists of three lines $l_0$, $l_1$ and $l_2$ which can be viewed as
coordinate axes for homogeneous coordinates in $\rp^2$.
Thus we have three affine charts for $\rp^2$. The intersection of all
three charts is $(\R^*)^2\subset\rp^2$. We denote $\R V=\R\bar{V}\cap(\R^*)^2$.
The complexification $V\subset (\C^*)^2$ is the complex hypersurface defined
by the same equation as $\R V$. Thus we are in position to apply the
content of the previous sections of the paper to the amoeba of $V$.

In \cite{Mi} it was shown (with the help of amoebas) that for each $d$ the topological type of the
pair $(\rp^2,\R\bar{V})$ is unique as long as the curve $\R\bar{V}$ is
maximal in each of the three affine charts of $\rp^2$. Furthermore,
the diffeomorphism type of the triad $(\rp^2;\R\bar{V},l_0\cup l_1\cup l_2)$
is unique. In subsection
\ref{krivye} we formulate this maximality condition and sketch the
proof of uniqueness. A similar statement holds for curves in other
toric surfaces. The Newton polygon $\Delta$ plays then the r\^ole of
the degree $d$.
In subsections \ref{poverxnosti} and \ref{hd} we describe an analogous
but weaker statement towards uniqueness.

\subsection{Relation to amoebas: the real part $\R V$ as a subset of the critical locus of $\Log|_V$
and the logarithmic Gauss map}

Suppose that the hypersurface $V\subset\tor$ is defined over real numbers (i.e. by a polynomial
with real coefficients).
Denote its real part via $\R V=V\cap\rtor$.
We also assume that $V$ is non-singular.
Let $F\subset V$ be the critical locus of the map $\Log|_V:V\to\R^n$.
It turns out that the real part $\R V$ is always contained in $F$.
\begin{prop}[Mikhalkin \cite{Mi}]
\label{RVF}
$\R V\subset F$.
\end{prop}
This proposition indicates that the amoeba must carry some information about $\R V$.
The proof of this proposition makes use of the {\em logarithmic Gauss map}.

Note that since $\tor$ is a Lie group there is a canonical trivialization
of its tangent bundle. If $z\in\tor$ then the multiplication by $z^{-1}$
induces an isomorphism $T_z\tor\approx T_1\tor$ of the tangent bundles at $z$
and $1=(1,\dots,1)\in\tor$.
\begin{defn}[Kapranov \cite{KaGauss}]
The {\em logarithmic Gauss map} is a map
$$\gamma:V\to\cp^{n-1}.$$ It sends each point $z\in V$ to the image
of the hyperplane $T_zV\subset T_z\tor$  under the canonical
isomorphism $T_z\tor\approx T_1\tor=\C^n$.

The map $\gamma$ is a composition of a branch of a holomorphic logarithm
$\tor\to\C^n$ defined locally up to translation by $2\pi i$ with the usual Gauss
map of the image of $V$.
We may define $\gamma$ explicitly in terms of the defining
polynomial $f$ for $V$ by logarithmic differentiation formula.
If $z=(z_1,\dots,z_n)\in V$ then
$$\gamma(z)=[<\nabla f,z>]=[\frac{\dd f}{\dd z_1}z_1:\dots:\frac{\dd f}{\dd z_n}z_n]\in\cp^{n-1}.$$
\end{defn}

\begin{lem}[Mikhalkin \cite{Mi}]
$F=\gamma^{-1}(\rp^{n-1})$
\end{lem}
To justify this lemma we recall that $\Log:\tor\to\R^n$ is a smooth fibration
and $V$ is non-singular.
Thus $z\in V$ is critical for $\Log|_V$ if and only if the tangent vector
space to $V$ and the tangent vector space to the fiber torus
$\gamma^{-1}(\gamma(z))$ intersect along an $(n-1)$-dimensional
subspace. Such points are mapped to real points of $\cp^{n-1}$ by $\gamma$.

Note that this lemma implies Proposition \ref{RVF}. If $V$ is defined over
$\R$ then $\gamma$ is equivariant with respect to the complex conjugation
and maps $\R V$ to $\rp^{n-1}$.

\subsection{Compactification: a toric variety associated to a hypersurface in $\tor$}
A hypersurface $V\subset\tor$ is defined by a polynomial $f:\C^n\to C$.
If the coefficients of $f$ are real then we define the real part of $V$
by $\R V=V\cap\rtor$. Recall that the Newton polyhedron $\Delta\subset\R^n$
of $V$ is an integer convex polyhedron obtained as
the convex hull of the indices of monomials participating in $f$,
see \eqref{NP} in subsection \ref{fpt}.

Let $\C T_\Delta\supset\tor$ be the toric variety corresponding to $\Delta$,
see e.g. \cite{GKZ} and let $\R T_\Delta\supset\rtor$ be its real part.
We define $\bar{V}\subset\C T_\Delta$ as the closure of $V$ in $\C T_\Delta$
and we denote via $\R\bar{V}$ its real part.

Note that $\bar{V}$ may be singular even if $V$ is not. Nevertheless $\C T_\Delta$
is, in some sense, the best toric compactification of $\tor$ for $V$. Namely, $\bar{V}$ does not
pass via the points of $\C T_\Delta$ corresponding to the vertices of $\Delta$ and
therefore it does not have singularities there. Furthermore, $\C T_\Delta$ is minimal
among such toric varieties, since $\bar{V}$ intersect any line in $\C T_\Delta$
corresponding to an edge of $\Delta$.

Thus we may naturally compactify the pair $(\tor,V)$ to the pair $(\C T_\Delta,\bar{V})$.
In such a setup the polyhedron $\Delta$ plays the r\^ole of the degree in $\C T_\Delta$.
Indeed, two integer polyhedra $\Delta$ define the same toric variety $\C T_\Delta$
if their corresponding faces are parallel. But the choice of $\Delta$ also fixes 
the homology class of $\bar{V}$ in $H_{2n-2}(\C T_\Delta)$.

The simplest example is the projective space $\cp^n$. The corresponding $\Delta$ is,
up to translation and the action of $SL_n(\Z)$ the simplex defined by equations
$z_j>0$, $z_1+\dots+z_n<d$. Thus in this case $\Delta$ is parametrized by a single
natural number $d$ which is the degree of $\bar{V}\subset\cp^n$.

\subsection{Maximality condition for $\R V$}
The inequality $l\le\frac{(d-1)(d-2)}{2}$ discovered by Harnack for the number $l$
of components of a curve $\R\bar{V}$ is a part of a more general {\em Harnack-Smith
inequality}.
Let $X$ be a topological space and let $Y$ be the fixed point set of a of a continuous
involution on $X$. Denote by $b_*(X;\Z_2)=\dim H_*(X;\Z_2)$ the total $\Z_2$-Betti number of $X$.
\begin{thm}[P. A. Smith, see e.g. the appendix in \cite{W}]
\label{smith}
$$b_*(Y;\Z_2)\le b_*(X;\Z_2).$$
\end{thm}
\begin{cor}
$b_*(\R\bar{V};\Z_2)\le b_*(\bar{V};\Z_2)$, $b_*(\R{V};\Z_2)\le b_*({V};\Z_2).$
\end{cor}
Note that Theorem \ref{smith} can also be applied to pairs which consist of
a real variety and real subvariety and other similar objects.
\begin{defn}[Rokhlin \cite{R}]
A variety $\R\bar{V}$ is called an {\em M-variety} if $$b_*(\R\bar{V};\Z_2)= b_*(\bar{V};\Z_2).$$
\end{defn}
E.g. if $\bar{V}\subset\cp^2$ is a smooth curve of degree $d$ then $\bar{V}$ is a Riemann surface
of genus $g=\frac{(d-1)(d-2)}{2}$.
Thus $b_*(\bar{V};\Z_2)=2+2g$. On the other hand, $b_*(\R\bar{V};\Z_2)=2l$, where $l$ is the number
of (circle) components of $\R\bar{V}$.

Let $\R V\subset\rtor$ be an algebraic hypersurface, $\Delta$ be its Newton polyhedron,
$\R T_{\Delta}$ be the toric variety corresponding to $\Delta$ and $\R\bar{V}\subset\R T_\Delta$
the closure of $\R V$ in $\R T_\Delta$.
We denote with $V\subset\tor$ and $\bar{V}\subset\C T_\Delta$ the complexifications of these
objects.
Recall (see e.g. \cite{GKZ}) that each (closed) $k$-dimensional face $\Delta'$ of $\Delta$ corresponds to 
a closed $k$-dimensional toric variety $\R T_{\Delta'}\subset\R T_\Delta$ (and, similarly,
$\C T_{\Delta'}\subset\C T_\Delta$). The intersection $V_{\Delta'}=\bar{V}\cap\C T_{\Delta'}$ is itself
a hypersurface in the $k$-dimensional toric variety $\C T_{\Delta'}$ with the Newton polyhedron $\Delta'$.
Its real part is $\R V_{\Delta'}=V_{\Delta'}\cap\R\bar{V}$.

Denote with $\operatorname{St}\Delta'\subset\dd\Delta$ the union of all the closed faces of $\Delta$
containing $\Delta'$. Denote
$V_{\operatorname{St}\Delta'}=\bigcup\limits_{\Delta''\subset\operatorname{St}\Delta'}V_{\Delta'}$
and $\R V_{\operatorname{St}\Delta'}=V_{\operatorname{St}\Delta'}\cap\R T_\Delta$.

\begin{defn}[Mikhalkin \cite{M-r}]
\label{max}
A hypersurface $\R\bar{V}\subset\C T_{\Delta}$ is called {\em torically maximal} if the following conditions
hold
\begin{itemize}
\item $\R\bar{V}$ is an M-variety, i.e. $b_*(\R\bar{V};\Z_2)=b_*(\bar{V};\Z_2)$;
\item the hypersurface $\bar{V}\cap\C T_{\Delta'}\subset\C T_{\Delta}$ is torically maximal for each
face $\Delta'\subset\Delta$ (inductively we assume that this notion is already defined in smaller
dimensions);
\item for each face $\Delta'\subset\Delta$ we have
$b_*(\R V\cup\R V_{\operatorname{St}\Delta'},\R V_{\operatorname{St}\Delta'};\Z_2)=
b_*(V\cup V_{\operatorname{St}\Delta'},V_{\operatorname{St}\Delta'};\Z_2)$.
\end{itemize}
\end{defn}

Consider a linear function $h:\R ^n\to\R$. A facet $\Delta'\subset\Delta$ is called {\em negative}
with respect to $h$ if the image of its outward normal vector under $h$ is negative.
We define
$\C T^-=\bigcup\limits_{\text{negative}\ \Delta'}\C T_{\Delta'}.$
In these formula we take the union over all the closed facets $\Delta'$ negative
with respect to $h$.
Let $V^-=\bar{V}\cap\C T^-$ and $\R V^-=V^-\cap\R\bar{V}$.

We call a linear function $h:\R^n\to\R$ generic if its kernel does not contain vectors orthogonal
to facets of $\Delta$.
\begin{prop}
If a hypersurface $\R\bar{V}\subset\R T_\Delta$ is torically maximal then for any generic linear function $h$
we have
$$b_*(\R V\cup\R V^-,\R V^-;\Z_2)=b_*(V\cup V^-,V^-;\Z_2).$$
\end{prop}

\subsection{Curves in the plane}
\label{krivye}
\subsubsection{Curves in $\rp^2$ and their bases}
Note that if $\R V\subset (\R^*)^2$ is a torically maximal curve then 
the number of components of $\R\bar{V}$ coincides with the genus of $\C \bar{V}$.
In other words (cf. \ref{H16}) $\R\bar{V}$ is an M-curve.

We start by reformulating the maximality condition of Definition \ref{max} for
the case of curves in the projective plane.
Let $\R C\subset\rp^2$ be a non-singular curve of degree $d$.

\begin{figure}[h]
\label{bbasa}
\centerline{
\psfig{figure=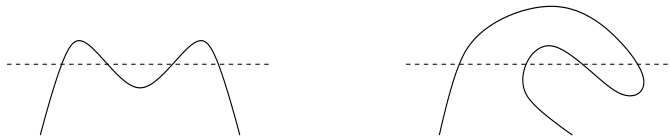,height=0.7in,width=3.5in}}
\caption{Possible bases for a real quartic curve.}
\end{figure}

\begin{defn}[Brusotti \cite{Br}]
Let $\alpha$ be an arc (i.e. an embedded closed interval) in $\R C$.
The arc $\alpha$ is called a {\em base} (or a {\em base of rank 1}, see \cite{Br})
if there exists a line $L\subset\rp^2$ such that the intersection $L\cap\alpha$
consists of $d$ distinct points.
\end{defn}

Note if three lines $L_1,L_2,L_3$ in $\rp^2$ are generic, i.e. they do not pass through
the same point, then =$\rp^2\setminus (L_1\cup L_2\cup L_3)=(\R^*)^2$. 
We call such $(\R^*)^2$ a {\em toric chart} of $\rp^2$.
Thus $\R V=\R C\setminus (L_1\cup L_2\cup L_3)$ is a curve in $(\R^*)^2$.
If $\R C$ does not pass via $L_j\cap L_k$ then the Newton polygon of $\R V$
(for any choice of coordinates $(x,y)$ in $(\R^*)^2$ extendable to affine coordinates
in $\R^2=\rp^2\setminus L_j$ for some $j$) is the triangle $\Delta_d=\{x\ge 0\}\cap
\{y\ge 0\}\cap\{x+y\le d\}$.
\begin{prop}[Mikhalkin \cite{M-r}]
The curve $\R C\subset\rp^2$ is maximal in some toric chart of $\rp^2$
if and only if $\R C$ is an M-curve with three disjoint bases.
\end{prop}
Many M-curves with one or two disjoint bases are known (see e.g. \cite{Br}).
However there is (topologically) only one known example of curve with three disjoint bases,
namely the first M-curve constructed by Harnack \cite{Ha}. Theorem \ref{urp2}
asserts that this example is the only possible.

\begin{defn}[simple Harnack curve in $\rp^2$, cf. \cite{Ha}, \cite{MR}]
\label{har}
A non-singular curve $\R C\subset\rp^2$ of degree $d$ is called a (smooth) simple Harnack curve if it is an M-curve and
\begin{itemize}
\item all ovals of $\R C$ are disjoint (i.e. have disjoint interiors, see \ref{H16}) if $d=2k-1$ is odd;
\item one oval of $\R C$ contains $\frac{(k-1)(k-2)}{2}$ ovals in its interior while all other ovals
are disjoint if $d=2k$ is even.
\end{itemize}
\end{defn}

\begin{figure}[h]
\label{deg10}
\centerline{
\psfig{figure=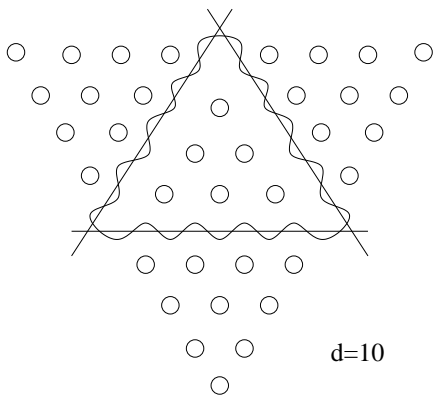,height=2.6in,width=2.6in}}
\caption{\cite{Mi} A simple Harnack curve.}
\end{figure}

\begin{thm}[Mikhalkin \cite{Mi}]
\label{urp2}
Any smooth M-curve $\R C\subset\rp^2$ with at least three base is a simple Harnack curve.
\end{thm}
There are several topological arrangements of M-curves with fewer than 3 bases 
for each $d$ (in fact, their number grows exponentially with $d$). There is a unique
(Harnack) topological arrangement of an M-curve with 3 bases by Theorem \ref{urp2}.
In the same time 3 is the highest number of bases an M-curve of sufficiently high degree
can have as the next theorem shows.
\begin{thm}[Mikhalkin \cite{Mi}]
No M-curve in $\rp^2$ can have more than 3 bases if $d\ge 3$.
\end{thm}

\subsubsection{Curves in real toric surfaces}
Theorem \ref{urp2} has a generalization applicable to other toric surfaces.
Let $\R V\subset(\R^*)^2$ be a curve with the Newton polygon $\Delta$. The sides of $\Delta$ correspond
to lines $L_1,\dots,L_n$ in $\R T_\Delta$. We have $\R V=\R\bar{V}\setminus(L_1\cup\dots\cup L_n$.
\begin{thm}[Mikhalkin \cite{Mi}]
\label{kriv}
The topological arrangement of a torically maximal curve is unique for each $\Delta$.
More precisely, the topological type of the triad $(\R T_\Delta;\R\bar{V},L_1\cup\dots\cup L_n)$
and, in particular, the topological type of the pair $((\R^*)^2,\R V)$ depends only on $\Delta$
as long as $\R V$ is a torically maximal curve.
\end{thm}

A torically maximal curve $\R\bar{V}$ is a counterpart of a simple Harnack curve for $\R T_\Delta$.
All of its components except for one are ovals with disjoint interiors. The remaining component
is not homologous to zero unless $\Delta$ is even (i.e. obtained from another lattice polygon
by a homothety with coefficient 2). If $\Delta$ is even the remaining component is also an oval
whose interior contains $g(V)$ ovals of $\R V$. Recall that, by Khovanskii's formula \cite{Kh},
$g(V)$ coincides with the number of lattice points in the interior of $\Delta$. 

\begin{thm}[Harnack, Itenberg-Viro \cite{Ha}, \cite{IV}]
For any $\Delta$ there exists a curve $\R V\subset(\R^*)^2$ which is torically maximal
and has $\Delta$ as its Newton polygon.
\end{thm}
As in Definition \ref{har} we call such curves {\em simple Harnack curves}, cf. \cite{MR}.

\subsubsection{Geometric properties of algebraic curves in $(\R^*)^2$}
It turns out that the simple Harnack curves have peculiar geometric properties,
but they are better seen after a logarithmic reparametrization $\Log|_{(\R^*)^2}:(\R^*)^2\to\R^2$.
A point of $\R V$ is called a logarithmic inflection point if it corresponds to an inflection
point of $\Log(\R V)\subset\R^2$ under $\Log$.

\begin{thm}[Mikhalkin \cite{Mi}]
\label{logkriv}
The following conditions are equivalent.
\begin{itemize}
\item $\R V\subset(\R^*)^2$ is a simple Harnack curve.
\item $\R V\subset(\R^*)^2$ has no real logarithmic inflection points.
\end{itemize}
\end{thm}

\begin{rem}[cf. \cite{Mi}]
Recall that by Proposition \ref{RVF} $\Log(\R V)$ is contained in the critical value
locus of $\Log|_V$. The map $\Log|_V:V\to\R^2$ is a surface-to-surface map in
our case and its most generic singularities are folds. By Proposition \ref{lc} the folds
are convex. Thus a logarithmic inflection point of $\R V$ must correspond to a
higher singularity of $\Log|_V$.

In \cite{Mi} it was noted that there are two types of stable (surviving small 
deformations of $\R V$) logarithmic inflection points of $\R V$. The first one,
{\em junction}, and correspond to intersection of $\R V$ with a
branch of imaginary folding curve. A junction logarithmic
inflection point can be found at the curve $y=(x-1)^2+1$. Note that the image of
the imaginary
folding curve under the complex conjugation is also a folding curve.
Thus over its image we have a double fold.

The second type, {\em pinching}, corresponds to intersection of $\R V$ with
a circle $E\subset V$ that gets contracted by $\Log$. Such circles $E$ survive
if we deform $V$ in the class of hypersurfaces with real coefficients but disappear
under a generic small perturbation if we allow the coefficients to become imaginary.

The circle $E$ intersect $\R V$ at two points. These points belong to different
quadrants of $(\R^*)^2$, but have the same absolute values of their coordinates.
Both of these points are logarithmic inflection points.

\begin{figure}[h]
\centerline{\psfig{figure=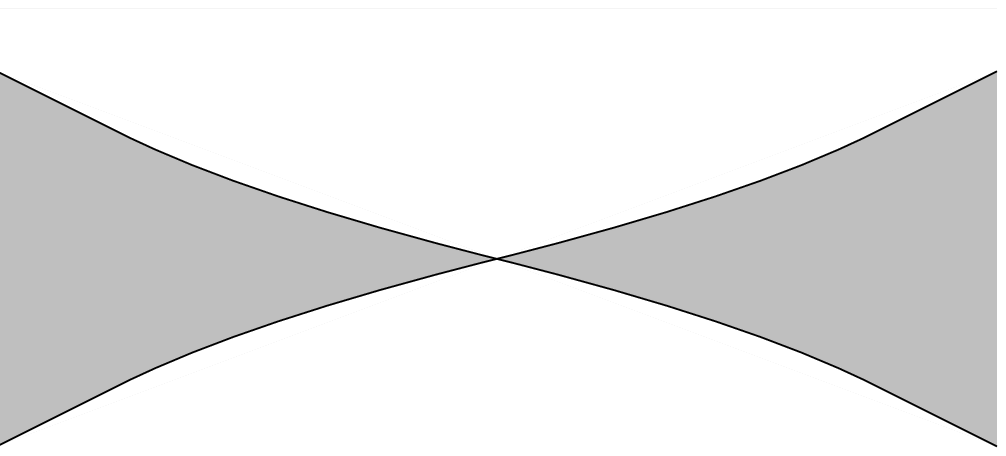,height=.8in,width=1.6in}\hspace{1in}
\psfig{figure=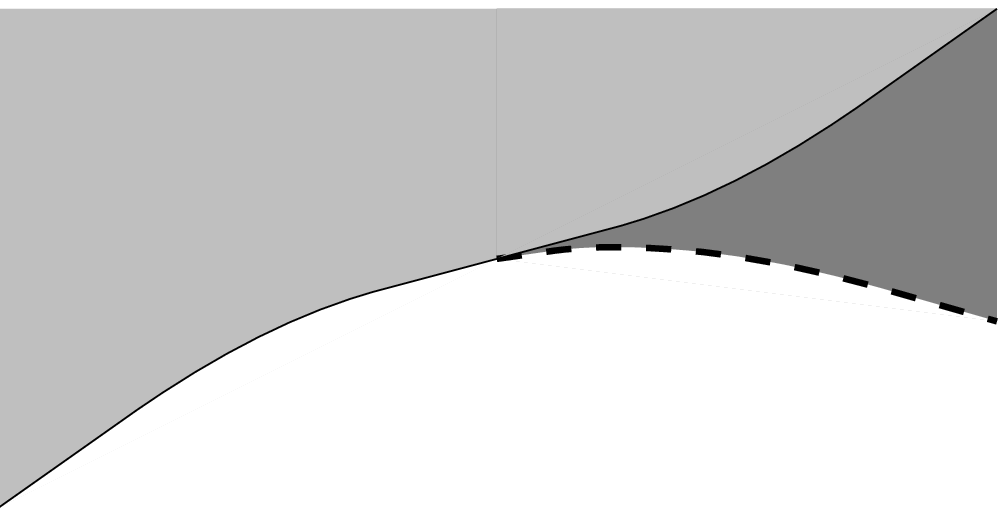,height=.6in,width=1.6in}}
\caption{\label{newsing} \cite{Mi} A pinching point and a junction point.}
\end{figure}
\end{rem}

\begin{prop} 
The logarithmic image $\Log(\R V)$ is trivial in the closed support homology
group $H^c_1(\R^2)$.
\end{prop}
Thus the curve $\Log(\R V)$ spans a surface in $(\R^*)^2$.
Theorem \ref{PR} has the following corollary.
\begin{cor}
The area of any region spanned by branches of $\Log(\R V)$ is smaller
than $\Area\Delta$.
\end{cor}

The situation is especially simple for
the logarithmic image of a simple Harnack curve.
\begin{prop}[\cite{Mi}]
\label{shemb}
If $\R V$ is a simple Harnack curve then $\Log|_{\R V}$ is an embedding
and $\Log\R V=\dd\am$.
\end{prop}

Thus in this case $\am$ coincides with the region spanned by the whole curve $\Log(\R V)$.
Furthermore, in \cite{MR} it was shown that simple Harnack curves maximize the
area of this region.
\begin{thm}[Mikhalkin-Rullg{\aa}rd, \cite{MR}]
If $\R V$ is a simple Harnack curve then
$\Area\am=\Area\Delta$.
\end{thm}

In the opposite direction we have the following theorem.
We say that a curve $V\subset(\C^*)^2$ is real up to translation if there exists
$a\in(\C^*)^2$ such that $aV$ is defined by a polynomial with real coefficients.
We denote the corresponding real part with $\R V$.
(Note that in general this real part might depend on the choice of translation.)
\begin{thm}[Mikhalkin-Rullg{\aa}rd \cite{MR}]
If $\Area\am=\Area\Delta>0$ and $V$ is non-singular and transverse to the
lines (coordinate axes) in $\C T_\Delta$ corresponding to the sides of $\Delta$
then $V$ is real up to translation in a unique way and $\R V$ is a simple Harnack curve.
\end{thm}
Furthermore, in \cite{MR} it was shown that the only singularities that $V$ can have
in the case $\Area\am=\Area\Delta>0$ are ordinary real isolated double points.

\subsection{Surfaces in the 3-space}
\label{poverxnosti}
\subsubsection{Topological uniqueness for torically maximal surfaces}
Let $\R V\subset(\R^*)^3$ be an algebraic surface with the Newton
polyhedron $\Delta\subset\R^3$. Let $\R\bar{V}\subset\R T_\Delta$
be its compactification.

Recall (see Definition \ref{max}) that if $\R V$ is a torically maximal surface then
$b_*(\R\bar{V};\Z_2)=b_*(\bar{V};\Z_2)$, i.e. $\R\bar{V}$ is an M-surface.

\begin{thm}[Mikhalkin \cite{M-r}]
\label{pov}
Given a Newton polyhedron $\Delta$ the topological type of
a torically maximal surface $\R\bar{V}\subset\R T_\Delta$ is unique.
\end{thm}
To describe the topological type of $\R\bar{V}$ it is useful to  compute
the total Betti number $b_*(\bar{V};\Z_2)$ in terms of $\Delta$. Note that
by the Lefschetz hyperplane theorem $b_*(\bar{V};\Z_2)=\chi(\bar{V})$.

We denote by $\Area\dd\Delta$ the total area of the faces of $\Delta$.
Each of these faces sits in a plane $P\subset\R^3$. The intersection  $P\cap\Z^3$
determines the area form on $P$. This area form is
translation invariant and such that the area of the smallest lattice
parallelogram is 1.

Similarly we denote by $\operatorname{Length}\operatorname{Sk}^1\Delta$
the total length of all the edges of $\Delta$. Again, each edge sits in a line $L\subset\R^3$.
The intersection $L\cap\Z^3$ determines the length on $L$ by setting the length of
the smallest lattice interval 1.

\begin{prop}
$b_*(V;\Z_2)=6\Vol\Delta-2\Area\dd\Delta+\operatorname{Length}\operatorname{Sk}^1\Delta.$
\end{prop}
This proposition follows from Khovanskii's formula \cite{Kh}.

\begin{thm}[Mikhalkin \cite{M-r}]
\label{ds}
A torically maximal surface $\R\bar{V}$ consists of $p_g+1$ components, where $p_g$
is the number of points in the interior of $\Delta$. There are $p_g$
components homeomorphic to 2-spheres and contained in $(\R^*)^3$.\
These spheres bound disjoint spheres in $(\R^*)^3$.
The remaining component is homeomorphic to
\begin{itemize}
\item a sphere with
$b_*(V;\Z_2)-2p_g(V)-2$ M\"obius bands in the case when $\Delta$ is odd
(i.e. cannot be presented as $2\Delta'$ for some lattice polyhedron $\Delta'$);
\item a sphere with $\frac12 b_*(V;\Z_2)-p_g(V)-1$ handles in the case $\Delta$ is even.
\end{itemize}
\end{thm}

\begin{rem}
\label{bertrand}
Not for every Newton polyhedron $\Delta$ a torically maximal surface $\R V\subset(\R^*)^3$
exists. The following example is due to B. Bertrand. Let $\Delta\subset\R^3$ be the convex hull
of $(1,0,0)$, $(0,1,0)$, $(1,1,0)$ and $(0,0,2k+1)$. If $k>0$ then there is no M-surface $\R\bar{V}$
with the Newton polyhedron $\Delta$. In particular, there is no torically maximal surface $\R V$
for $\Delta$.
\end{rem}

\begin{exa}
There are 3 different topological types of smooth M-quartics in $\rp^3$ (see \cite{Kha}).
They realize all topological possibilities for maximal real structures on abstract
K3-surfaces. Namely, such real surface may be homeomorphic to
\begin{itemize}
\item the disjoint union of 9 spheres and a surface of genus 2;
\item the disjoint union of 5 spheres and a surface of genus 6;
\item the disjoint union of a sphere and a surface of genus 10.
\end{itemize}
Theorem \ref{ds} asserts that only the last type can be a torically maximal quartic in $\rp^3$.
More generally, only the last type can be a torically maximal surface is a toric 3-fold
$\R T_\Delta$.
\end{exa}

\subsubsection{Geometric properties of maximal algebraic surfaces in $(\R^*)^3$}
Recall the classical geometric terminology.
Let $S\subset\R^3$ be a smooth surface.
We call a point $x\in S$ {\em elliptic}, {\em hyperbolic} or {\em parabolic}
if the Gauss curvature of $S$ at $x$ is positive, negative or zero.
\begin{rem}
Of course we do not actually need to use the Riemannian metric on $S$
do define these points. Here is an equivalent definition without referring
to the curvature. Locally near $x$ we can present $S$ as the graph of a function
$\R^2\to\R$. If the Hessian form of this function at $x$ is degenerate then we
call $x$ parabolic. If not, the intersection of $S$ with the tangent plane at $x$
is a real curve with an ordinary double point in $x$. If this point is isolated
we call $x$ elliptic. If it is an intersection of two real branches of the curve
we call it hyperbolic.
\end{rem}

We say that a point $x\in\R V\subset (\R^*)^3$ is {\em logarithmically} elliptic,
hyperbolic or parabolic if it maps to such point under
$\Log|_{(\R^*)^3}:(\R^*)^3\to\R^3$.

Generically for a smooth surface in $\R^3$ the parabolic locus, i.e. the set of parabolic points,
is a 1-dimensional curve. So is the logarithmic parabolic locus for a surface in $(\R^*)^3$.
In a contrast to this we have the following theorem for torically maximal surfaces.
Note that torically maximal surfaces form an open subset in the space of all surfaces
with a given Newton polyhedron.

\begin{thm}[Mikhalkin \cite{M-r}]
\label{parpov}
The logarithmic parabolic locus of a torically maximal surface
consists of a finite number of points.
\end{thm}
Note that such a zero-dimensional locus cannot separate the surface $\R V$.
Thus each component of $\R V$ is either logarithmically elliptic (all its points except finitely
many are logarithmically elliptic) or logarithmically hyperbolic
(all its points except finitely many are logarithmically hyperbolic).
\begin{cor}[Mikhalkin \cite{M-r}]
Every compact component of $\R V$ is diffeomorphic to a sphere.
\end{cor}
This corollary is a part of Theorem \ref{ds}.

\begin{rem}[logarithmic monkey saddles of $\R V$]
The Hessian at the isolated parabolic points $\Log(\R V)$ vanishes.
Generic parabolic points sitting on hyperbolic components of $\Log(\R V)$
look like so-called Monkey saddles (given in some local coordinates $(x,y,z)$
by $z=x(y^2-x^2)$).

Logarithmic monkey saddles do not appear on generic {\em smooth} surfaces in $(\R^*)^3$.
But they do appear on generic {\em real algebraic} surfaces in $(\R^*)^3$.
In particular, they appear on every torically maximal surface of sufficiently high degree.

The counterpart on the elliptic components of $\Log(\R V)$, the {\em imaginary monkey saddles},
are locally given by $z=x(y^2+x^2)$.
\end{rem}

\subsection{Hypersurfaces of higher dimension}
\label{hd}
Let $\R V\subset\rtor$ be a hypersurface and $n\ge 4$.
Theorems \ref{kriv} and \ref{pov} have a weaker version for
these dimensions.

\begin{thm}[Mikhalkin \cite{M-r}]
If $\R V$ is torically maximal then every compact component of $\R V$ is a sphere.
All these $(n-1)$-spheres bound disjoint $n$-balls in $\rtor$.
\end{thm}

The following theorem is a counterpart of Theorem \ref{parpov} and a weaker version
of Theorem \ref{logkriv}.
\begin{thm}[Mikhalkin \cite{M-r}]
The parabolic locus of $\Log(\R V)\subset\R^n$ is of codimension 2 if $\R V$ is torically maximal.
\end{thm}

\subsection{Maximality conditions for non-Archimedian amoebas}
Let $\am_K\subset\R^n$ be a non-Archimedian amoeba (see \ref{nA})
whose Newton polyhedron is $\Delta$.
\begin{prop}
The number of vertices of $\am_K$ is not greater than $n!\Vol\Delta$.
\end{prop}
This proposition can be deduced from Theorem \ref{thmka} and the fact
that the smallest possible volume of a convex lattice polyhedron is $\frac{1}{n!}$.

\begin{defn}
\label{nAmax}
A non-Archimedian amoeba $\am_K$ is called {\em maximal} if the number of its
vertices equals to $n!\Vol\Delta$.
\end{defn}

\begin{rem}
\label{remnAmax}
For some choices of $\Delta$ maximal amoebas do not exist.
We can take, for instance, $\Delta\subset\R^3$ to be the tetrahedron
with vertices $(1,0,0)$, $(0,1,0)$, $(0,0,1)$ and $(0,0,k)$. Any valuation
on the corresponding coefficients would have to be linear. Its Legendre
transform would have just one vertex while $\Vol\Delta=\frac{k}{n!}$.
Note that these polyhedra were used by B. Bertrand to show
absence of real maximal surfaces, see Remark \ref{bertrand}.
 
Nevertheless, if the toric variety corresponding to $\Delta$ is a projective space
or a product of projective spaces then maximal non-Archimedian amoebas exist.
This statement is implicitly contained in \cite{IV}.
\end{rem}

\section{Patchworking of amoebas, Maslov's dequantization
and topology of complex algebraic varieties}
\label{secmas}
\subsection{Patchworking polynomial}
\label{pp}
In 1979 Viro discovered a {\em patchworking} technique for
construction of real algebraic hypersurfaces, see \cite{V}.

Fix a convex lattice polyhedron $\Delta\in\R^n$.
Choose a function $v:\Delta\cap\Z^n\to\R$.
The graph of $v$ is a discrete set of points in $\R^n\times\R$.
The overgraph is a family of parallel rays. Thus the convex hull
of the overgraph is a semi-infinite polyhedron $\tilde\Delta$.
The facets of $\tilde\Delta$ which project isomorphically to
$\R^n$ define a subdivision of $\Delta$ into smaller
convex lattice polyhedra $\Delta_k$.

Let $F(z)=\sum\limits_{j\in\Delta}a_jz^j$ be a generic polynomial in the class of
polynomial whose Newton polyhedron is $\Delta$.
The {\em truncation} of $F$ to $\Delta_k$ is $F_{\Delta_k}=\sum\limits_{j\in\Delta_k}a_jz^j$.
The {\em patchworking polynomial} $f$ is defined by formula
\begin{equation}
\label{patch}
f^v_t(z)=\sum\limits_j a_jt^{v(j)}z^j,
\end{equation}
$z\in\R^n$, $t>1$ and $j\in\Z^n$.

Consider the hypersurfaces $V_{\Delta_k}$ and $V_t$ in $\tor$ defined
by $F_{\Delta_k}$ and $f^v_t$.
If $F$ has real coefficients then we denote $\R V_{\Delta_k}=V_{\Delta_k}\cap\rtor$
and $\R V_t=V_t\cap\rtor$. 
Viro's patchworking theorem \cite{V} asserts that for large values of $t$
the hypersurface $\R V_t$ can be obtained from $\R V_{\Delta_k}$ by a certain patchworking
procedure. The same holds for amoebas of the hypersurfaces $V_t$ and $\R V_{\Delta_k}$.
In fact patchworking of real hypersurfaces can be interpreted as the real version
of patchworking of amoebas (cf. Appendix in \cite{Mi}).
Below we describe a special case of amoeba's patchworkings
in terms of the so-called dequantization.

\subsection{Maslov's dequantization}
\label{deq}
It was noted by Viro in \cite{V2000} that patchworking is related to so-called
{\em Maslov's dequantization} of positive real numbers.

Recall that a {\em quantization} of a semiring $R$ is a family of semirings $R_h$,
$h\ge 0$ such that $R_0=R$ and $R_t\approx R_s$ as long as $s,t>0$, but
$R_0$ is not isomorphic to $R_t$.
The semiring $R_h$ with $h>0$ is called a {\em quantized} version of $R_0$.

Maslov (see \cite{Ma}) observed that the ``classical" semiring $\R_+$
of real positive number is a quantized version of some other ring in this sense
Let $R_h$ be the set of positive numbers with the usual multiplication
and with the addition operation 
$z\oplus_h w=(z^{\frac{1}{h}}+w^{\frac{1}{h}})^h$
for $h>0$ and
$z\oplus_h w=\max\{z,w\}$
for $h=0$.
Note that
$$\lim\limits_{h\to 0} (z^{\frac{1}{h}}+w^{\frac{1}{h}})^h=\max\{z,w\}$$
and thus this is a continuous family of arithmetic operations.

The semiring $R_1$ coincides with the standard semiring $\R_+$.
The isomorphism between $\R_+$ and $R_h$ with $h>0$ is given
by $z\mapsto z^h$. On the other hand the semiring $R_0$ is not
isomorphic to $\R_+$ since it is idempotent, indeed $z+z=\max\{z,z\}=z$.

\subsection{Logarithmic dequantization}
\label{logdeq}
Alternatively we may define the dequantization deformation with the help
of the logarithm.
The logarithm $\log_t$, $t>1$, induces a semiring structure on $\R$ from $\R_+$,
$$x\oplus_t y=\log_t(t^x+t^y),\ x\otimes_t y=x+y,\ x,y\in\R.$$
Similarly we have $x\oplus_\infty y=\max\{x,y\}$.
Let $R^{\log}_t$ be the resulting semiring.
\begin{prop}
The map $\log:R_h\to R^{\log}_t$, where $t=e^{\frac1h}$, is an isomorphism.
\end{prop}

\subsection{Patchworking as a dequantization}
The patchworking polynomial \eqref{patch} can be viewed as
a deformation of the polynomial $f^v_1$. We define a similar
deformation with the help of Maslov's dequantization.
Instead of deforming the coefficients we keep coefficients the
same and deform arithmetic operations as in \ref{deq} and \ref{logdeq}.

Choose any coefficients $\alpha_j$, $j\in\Delta$.
Let $\phi_t:(R^{\log}_t)^n\to R^{\log}_t$, $t\ge e$, be a polynomial whose
coefficients are $\alpha$, i.e. $$\phi_t(x)=\bigoplus_t(\alpha_j+jx),\ x\in\R^n.$$
Let $\Log_t:\tor\to\R^n$ be defined by $(x_1,\dots,x_n)=(\log|z_1|,\dots,\log|z_n|)$.
\begin{prop}[Maslov \cite{Ma},Viro \cite{V2000}]
The function $f_t=(\log_t)^{-1}\circ\phi_t\circ\Log_t:(\R_+)^n\to\R_+$ is a polynomial with respect to the standard
arithmetic operations in $\R_+$,
$$f_t(z)=\sum\limits_j t^{\alpha_j} z^j.$$
\end{prop}

This is a special case of the patchworking polynomial \eqref{patch}. The coefficients
$\alpha_j$ define the function $v:\Delta\cap\Z^n\to\R$.

\subsection{Limit set of amoebas}
Let $V_t\subset\tor$ be the zero set of $f_t$ and let $\am_t=\Log_t(V_t)\subset\R^n$.
Note that $\am_t$ is the amoeba of $V_t$ scaled $\log t$ times.
Note also that the family $f_t=\sum\limits_j t^{\alpha_j} z^j$ can be considered
as a single polynomial whose coefficients are powers of $t$. In particular
we may treat it as a polynomial over the field of Puiseux series, i.e. a
non-Archimedian field (see Example \ref{pui}). Let $\am_K$ be the corresponding
non-Archimedian amoeba.

We have a uniform convergence of the addition operation
in $R^{\log}_t$ to the addition operation in $R^{\log}_\infty$.
As it was observed by Viro it follows from the following inequality
$$\max\{x,y\}\le x\oplus_t y=\log_t(t^x+t^y)\le \max\{x,y\}+\log_t 2.$$
More generally, we have the following lemma.
\begin{lem}
$$\max\limits_{j\in\Delta}(\alpha_j+jx) \le \phi_t(x)\le \max\limits_{j\in\Delta}(\alpha_j+jx)+
\log N,$$
where $N$ is the number of lattice points in $\Delta$.
\end{lem}
The following theorem is a corollary of this lemma.

\begin{thm}[Mikhalkin \cite{M-l}, Rullg{\aa}rd \cite{R2}]
\label{amlim}
The subsets $\am_t\subset\R^n$ tend in the Hausdorff metric to $\am_K$.
\end{thm}

Note that by Theorem \ref{thmka} $\am_K$ is obtained by patchworking of
the amoebas of the truncations of $f_e$ to smaller polyhedra $\Delta_k$ (see \ref{pp}).

\subsection{Torus fibrations for algebraic hypersurfaces in $\tor$}
Recall that as long as a hypersurface $V\cap\tor$ is non-singular
its diffeomorphism type depends only on its Newton polyhedron $\Delta$.
Theorem \ref{amlim} implies that for large values of $t$ the amoeba
$\am_t$ is contained in a regular neighborhood $W$ of $\am_K$.

The space $\am_K$ has a natural cellular decomposition which turns $\am_K$ to
an $(n-1)$-dimensional CW-complex. The decomposition
comes from piecewise-linear embedding of $\am_K$ into $\R^n$ (cf. Theorem \ref{thmka}).
Each $k$-cell of $\am_K$ is contained in an affine $k$-subspace of $\R^n$.

\begin{prop}[Mikhalkin \cite{M-l}]
Let $z\in\am_K$ be a point of an open $(n-1)$-cell.
Let $\rho:W\to\am_K$ be a regular neighborhood retraction such that
its restriction to a neighborhood of $z$ is a smooth submersion.
For sufficiently large $t>0$
the composition $$\lambda:V_t\stackrel{\Log_t}\to W\stackrel\rho\to\am_K$$
is submersive near $z$ and $\lambda^{-1}(z)$ is diffeomorphic to 
a smooth $(n-1)$-torus.
\end{prop}

If $\am_K$ is maximal (see Definition \ref{nAmax}) then
the map $\lambda$ can be further improved. Let $z\in\am_K$.
\begin{defn}
\label{degl}
Let $M$ be a manifold, $N\subset\R^n$ be a piecewise-smooth CW-complex
and $\lambda:M\to N\subset\R^n$ be a smooth map.
Let $x\in L$ be a point.
By the {\em degeneration type} of $\lambda$ near $x$ we mean the equivalence class
of the restriction
$\lambda^{-1}(U)\to U$ of $\lambda$ to a small open ball $U=N\cap D_x(\epsilon)$
near $x$. Two smooth maps $W\to U\subset D_x(\epsilon)$ and
$W'\to U'\subset D_{x'}(\epsilon')$ are equivalent if there exist diffeomorphisms
$W\stackrel{\approx}{\to} W'$ and $D_x(\epsilon)\stackrel{\approx}{\to} D_{x'}(\epsilon')$
which take the first map to the second map.
\end{defn}


\begin{thm}[Mikhalkin \cite{M-l}]
\label{tortor}
Suppose that $\am_K$ is maximal.
There exists a regular neighborhood retraction $\rho:W\to\am_K$
such that for sufficiently large $t>0$ the composition
$$\lambda=\rho\circ(\Log_t|_{V_t}):V_t\to\am_K$$
is a singular torus fibration in the following sense
\begin{itemize}
\item the restriction of $\lambda$ to any open cell of $\am_K$ 
is a trivial fibration;
\item the fiber of $\lambda$ over an $(n-1)$-cell is $T^{n-1}$;
\item the degeneration type of $\lambda$ at $x\in\am_K$
depends only on the dimension of the open cell containing $x$.
\end{itemize}
The fiber of $\lambda$ over an open $k$-cell of $\am_K$ is
a $(n-1)$-dimensional CW-complex that can be embedded to
the $n$-torus $T^n$. The fiber over an open $(n-2)$-cell
is the product of a $\theta$-graph (i.e. the graph with 2 vertices
and 3 edges joining them) and a torus $T^{n-2}$.
In addition we have the following properties.
\begin{itemize}
\item The base $\am_K$ is homotopy equivalent
 to a wedge of $p_g$ spheres $S^{n-1}$, where $p_g$ is the number
of lattice points in the interior of $\Delta$. 
\item The induced homomorphism
$$\lambda^*:H^{n-1}(\am_K;\Z)\to H^{n-1}(V;\Z)$$
is a monomorphism.
\end{itemize}
\end{thm}

\subsection{Torus fibrations for complex projective hypersurfaces}
This theorem admits a compactified version.
Let $\bar{V}\subset\C T_\Delta$ be the compactification of $V$.
A non-Archimedian amoeba corresponding to $\Delta$ can be
compactified as well. Recall (see Remark \ref{repar})
that the moment maps for the symplectic spaces $\tor$ and $\C T_\Delta$ 
define a reparametrization $\R^n\stackrel{\approx}{\to}\Int\Delta$.
The {\em compactified non-Archimedian amoeba} $\Pi$ is the closure in $\Delta$
of the image of a non-Archimedian amoeba under this reparametrization.

Note that $\Pi$ admits a natural cellular structure. To each cell we can associate
two indices. One index is its dimension $k$. The other is the dimension $l$ of the (open)
face of $\Delta$ containing the cell. 
 
\begin{defn}
\label{specspine}
An $(n-1)$-dimensional cellular space $\Pi$ is called a {\em special spine}
if a small neighborhood of a point
$x\in\Pi$ from an open $k$-cell is homeomorphic to the direct product
of $\R^k$ and the cone over the $(n-k-2)$-skeleton of the $(n-k)$-dimensional
simplex.

The space $\Pi$ is called a {\em special spine with corners} if for each open
$k$-dimensional cell there exists an integer number $l$, $k<l\le n$ with the
following property. A small neighborhood of a point $x\in\Pi$ from this cell
is homeomorphic to the direct product of $\R^k\times [0,+\infty)^{n-l}$ and
the cone over the $(l-k-2)$-skeleton of the $(l-k)$-dimensional
simplex. Note that a $(-1)$-skeleton is always empty.
Such a $k$-dimensional cell is called a {\em $(k,l)$-cell}.
\end{defn}
\begin{exa}
A 1-dimensional special spine is a 3-valent graph.
A 1-dimensional special spine with corners is a 3- and 1-valent graph.
\end{exa}

\begin{prop}
If $\am_K$ is a maximal non-Archimedian amoeba then $\Pi$ is a special spine
with corners.
\end{prop}
\begin{rem}
The term ``special spine" comes from Topology. 
Let $X$ be an $n$-manifold (possibly with boundary or even with corners).
An $(n-1)$-dimensional CW-complex $S\subset X$
is called a spine of $X$ if the complement $X\setminus (S\cup\dd X)$
is a disjoint union of open $n$-balls.

Originally the term ``special spine" referred to a spine which satisfies
to additional properties specified in Definition \ref{specspine}.
Now the this term is also used (in particular, in this paper) also for
CW-complexes without any ambient space. Note that in our case
$\Pi$ is a spine of $\Delta$ in the topological sense.
\end{rem}

We introduce the following definition for the next theorem.
\begin{defn}[Mikhalkin \cite{M-l}]
A map $\lambda:M\to\Pi$ is called a {\em manifold fibration over a
special spine $\Pi\subset\Delta$ with corners}, where $\Delta\subset\R^n$
is a convex polyhedron, it
\begin{itemize}
\item $M$ is a manifold;
\item $\lambda:M\to\Pi\subset\Delta$ is a smooth map;
\item the restriction of $\lambda$ to each open cell of $\Pi$ is
a smooth trivial fibration (a submersion over an open cell);
\item the degeneration type (see Definition \ref{degl}) of $\lambda$ at a point $x$
from an open $(k,l)$-cell depends only on $k$ and $l$.
\end{itemize}
\end{defn}
Note that an $(n-1)$-dimensional cell is always a $(n-1,n)$-cell.

Remark \ref{remnAmax} states that maximal non-Archimedian
amoebas exist in the case when $\C T_\Delta$ is a projective space
or a product of projective spaces.
If $\am_K$ is maximal then the
corresponding compactified non-Archimedian amoeba $\Pi$ is

\begin{thm}[Mikhalkin \cite{M-l}]
\label{torcomp}
Let $\bar{V}\subset\cp^n$ be a non-singular hypersurface.
There exists a special spine $\Pi$ with corners and a
manifold fibration $\bar\lambda:\bar{V}\to\Pi$ over $\Pi$ such that
\begin{itemize}
\item the general fiber of $\bar\lambda$ (i.e. the fiber over an open $(n-1)$-dimensional cell)
is a smooth $(n-1)$-dimensional torus;
\item the homotopy type of $\Pi$ is the wedge of $p_g$ copies of $S^{n-1}$,
where $p_g=h^{n-1,0}$ is the geometric genus of $\bar{V}$;
\item the induced homomorphism
$\lambda^*:H^{n-1}(\Pi;\Z)\to H^{n-1}(\bar{V};\Z)$ is a monomorphism.
\end{itemize}
\end{thm}

\begin{add}[Mikhalkin \cite{M-l}]
Here is a partial description of special fibers of $\bar\lambda$ from Theorem \ref{torcomp}.
\begin{itemize}
\item The fiber of $\bar\lambda$ over an $(k,k+1)$-cell, $k<n$, is a smooth $k$-dimensional torus;
\item the fiber of $\bar\lambda$ over an $(k,k+2)$ cell, $k<n-1$, is a product of the $\theta$-graph
(i.e. the graph with 2 vertices and 3 edges joining them) and a $(k-1)$-torus;
\item more generally, the fiber of $\bar\lambda$ over a $(k,l)$-cell is an $(l-1)$-dimensional CW-complex
whose topology depends only on $k$ and $l$ and such that it can be
embedded to an $l$-dimensional torus.
\end{itemize}
\end{add}

\begin{add}[Mikhalkin \cite{M-l}]
\label{addpants}
Let $x\in\Pi$ be a point from a $(k,l)$-cell and $U\ni x$ be a regular neighborhood of $x$ in $\Pi$.
The inverse image $\bar\lambda^{-1}(U)$ is diffeomorphic to the product of
$\R^k\times [0,+\infty)^{n-l}$ and $\C P^{l-k-1}$ minus
$l-k+1$ hyperplanes in general position. 
\end{add}

\subsection{Decomposition of projective hypersurfaces into pairs of pants}
Let $S$ be a closed Riemann surface.
An {\em open pair of pants} is an open manifold diffeomorphic to the two sphere $S^2$ 
minus 3 points.
A (closed) {\em pair of pants} is a compact surface of genus 0 with 3 boundary components.
It is easy to see that an open pair of pants is a pair of pants without its boundary.

A {\em pair of pants decomposition for} $S$ is given by a collection 
of disjoint embedded circles such that each component of their complement
is an open pair of pants.

Let $p_1,\dots,p_m\in S$ are distinct points.
A {\em pair of pants decomposition for} $(S;p_1,\dots,p_m)$ is given by a collection 
of disjoint embedded circles such that each component of their complement
in $S\setminus\bigcup\limits_j\{p_j\}$ is an open pair of pants.

\begin{prop}
To a pair of pants decomposition of $S$ we may canonically associate a
manifold fibration $\lambda:S\to\Pi$ over a 3-valent graph $\Pi$.

To a pair of pants decomposition of $(S;p_1,\dots,p_m)$ we may canonically associate a
manifold fibration $\lambda:S\to\Pi$ over a 3- and 1-valent graph $\Pi$.
\end{prop}
Note that there is a natural fibration of a pair of pants over
a Y-shaped graph such that the boundary components are fibers over
1-valent vertices and the fiber over the 3-valent vertex is a $\theta$-shaped graph.

In the opposite direction we have the following proposition.
\begin{prop}
Let $S\to\Pi$ be a manifold fibration over a 3-valent graph $\Pi$ such that
the fibers over 3-valent vertices are $\theta$-shaped graphs. Then
the inverse images of the midpoints of the edges give a pair of pants decomposition
for $S$.

Let $S\to\Pi$ be a manifold fibration over a 3- and 1-valent graph $\Pi$ such that
the fibers over 3-valent vertices are $\theta$-shaped graphs and the fibers over
1-valent vertices are points $p_1,\dots,p_m$. Then
the inverse images of the midpoints of the edges connecting 3-valent vertices
give a pair of pants decomposition for $(S;p_1,\dots,p_m)$.
\end{prop}
The graph $\Pi$ can be interpreted as combinatorial data needed for gluing pairs of
pants to obtain $S$.
\begin{prop}
\label{1dimreconstruct}
The surface $S$ may be recovered from $\Pi$ by the following procedure.
\begin{enumerate}
\item Take a disjoint union of pairs of pants, one pair of pants for each 3-valent vertex
of $\Pi$.
\item For each edge connecting 3-valent vertices identify some boundary components
of the corresponding pairs of pants. 
\item Collapse the remaining boundary components (those corresponding to
1-valent vertices to points.
\end{enumerate}
\end{prop}

\begin{defn}[Mikhalkin \cite{M-l}]
An open $l$-dimensional pair of pants is an open manifold diffeomorphic to
$\cp^l$ minus $l+2$ hyperplanes in general position.
\end{defn}
Note that the arrangement of $l+2$ hyperplanes in general position is
unique up to the natural action of $PSL(l+1,\C)$.
\begin{defn}[Mikhalkin \cite{M-l}]
An $l$-dimensional pair of pants $P_l$ is an compact manifold (with corners)
diffeomorphic to $\cp^l$ minus the union of small tubular neighborhoods $l+2$
hyperplanes in general position.
A {\em closed facet} of $P_l$ is the intersection of $\dd P_l$ and the boundary
of the tubular neighborhood of one of the $l+2$ hyperplanes.
A {\em closed $m$-face} of $P_l$ is the intersection of $l-m$ facets in $\dd P_l$.
An {\em open $m$-face} is a closed $m$-face minus all smaller-dimensional faces.
\end{defn}
Note that an open $m$-face of $P_l$ is an open manifold diffeomorphic to 
the open $m$-dimensional pair of pants $P_m$ times the real $(l-m)$-torus $T^{l-m}$.
Note also that an open pair of pants is a pair of pants minus its boundary.
\begin{rem}
\label{dcollapse}
We can collapse a part of the boundary of $P_l$ corresponding to a $m$
facets of $P_l$. The result of collapse is $\cp^l$ minus the union of small
tubular neighborhoods of the remaining $l+2-m$ hyperplanes. Thus we
add back the tubular neighborhoods of the hyperplane corresponding
to collapsing facets. 
\end{rem}

Theorem \ref{torcomp} can be interpreted as a higher-dimensional pair of pants
decomposition for smooth projective hypersurfaces thanks to the following
corollary from Addendum \ref{addpants}.
\begin{cor}
\label{pants}
Let $x\in\Pi$ be a $(0,n)$-cell of $\Pi$
and $U\ni x$ be a regular neighborhood of $x$.
The inverse image $\bar\lambda^{-1}(U)$ is diffeomorphic to
an open $(n-1)$-dimensional pair of pants.
\end{cor}

The polyhedral complex $\Pi$ may be interpreted as combinatorial data
needed for gluing pairs of pants to construct $\bar{V}$ in a fashion similar
to Proposition \ref{1dimreconstruct}.
We start from a disjoint union of $(n-1)$-dimensional pairs of pants, one
for each $(0,n)$-cell of $\Pi$. Each $(1,n)$-cell is an edge connecting
$(0,n)$-vertices.

Each $(0,n)$-vertex is adjacent to $n+1$ edges of $\Pi$ corresponding to $n+1$
facets of $P_{n-1}$. Similarly, it is adjacent to
$\begin{pmatrix} n+1 \\ k \end{pmatrix}$ $k$-faces of $\Pi$ corresponding to
$\begin{pmatrix} n+1 \\ k \end{pmatrix}$ $(n-k-1)$-faces of $P_{n-1}$.
For each $(1,n)$-edge we identify corresponding closed facets of the
pairs of pants corresponding to the endpoints.

Our identification is subject to the following additional condition.
For each $k$-cell $e$ of $\Pi$ we consider all $(0,n)$-vertices adjacent
to $e$. Each of the corresponding $(n-1)$-dimensional pair of pants
contain an $(n-k-1)$-faces corresponding to $e$. All these $(n-k-1)$-faces
have to be identified.

To get $\bar{V}$ from the result of this identification we have to collapse
the boundary as in Remark \ref{dcollapse}.

\end{document}